\documentclass[a4paper,11pt]{article}

\usepackage[latin1]{inputenc}   
\usepackage[T1]{fontenc} 
\usepackage{amsmath}
\usepackage[dvips]{graphicx}
\usepackage{theorem}
\usepackage{amsfonts}

\newcommand\om{\omega}
\newcommand\x{\chi}

\newtheorem{thm}{Theorem}[section]
\newtheorem{prop}{Proposition}[section]

\title {\textbf{Distinguished dihedral representations of $GL(2)$ over a p-adic field}}

\begin{document}
\maketitle

\section{Introduction}

Let $F$ be a finite extension of a p-adic field, $K$ a quadratic extension of $F$.\\
The principal series representations of $GL _2 (K)$ distingushed for $GL _2(F)$ are well known, it's also known that the Steinberg representation is distingushed ( cf.\cite{AT} for a summary of these results due to Y.Flicker, J.Hakim and D.Prasad,).\\
Moreover there is a complete characterisation of distinguishedness in terms of the epsilon factor (due to J.Hakim) and in terms of base change from representations of the unitary group U(2, K/F) (due to Y.Flicker, cf.\cite{AT} for a local proof).\\
Using those, we give here a description of dihedral representations on the parameter's side (those which the Langlands correpondance associate with 2 dimension induced representations of the Weil group of $K$, cf.(2) p.122 in \cite{GL}) distinguished with respect to $GL _2(F)$. \\
Every distinguished representation of $GL _2 (K)$ is paramatrised by a regular multiplicative character $\omega$ of a quadratic extension $L$ of $K$.\\
 We show ( theorem 5.1) that such a representation is distinguished for $GL _2(F)$ if and only if one can choose $L$ to be biquadratic over $F$ and $\omega$ trivial one the invertible elements of the two other quadratic extensions of $F$ in $L$.\\
The results we prove here are theorems 3.3, 5.1 and 6.1.
The method is to isolate first representations distingushed for ${GL _2(F)} _{+}$ using theorem 3.1, then to determinate those who are $GL _2(F)$-distinguished using theorems 4.1 and proposition 4.1.\\
Thus, in the case of odd residual caracteristic, we obtain every supercuspidal distinguished representations.\\
We also observe ( see proposition 5.3) that if we consider the principal series as paramatrised by a multiplicative character of a two dimensional semi-simple commutative algebra over $K$, the statement for distinguishedness is the same as for the supercuspidal dihedral representations.\\
We then give a generalisation of theorem 5.1 to dihedral representations ( non necessarily supercuspidal) in theorem 6.1.

\section{Preliminaries}

\subsection{Generalities}
 
 We consider $F$ a finite extension of ${\mathbb{Q}}_p$, and $K$ a quadratic extension of $F$ in an algebraic closure $\bar{F}$ of ${\mathbb{Q}}_p$.\\
 If $L$ is a quadratic extension of $K$ in $\bar{F}$, then to every character $\omega$ of $L^*$, we associate a representation of $GL_{2}(K)$ via the Weil representation ( cf.\cite{JL} p.144). \\
Such a representation is called dihedral. \\
We note $\theta$ the conjugation of $F/K$. \\
For $A$ is a ring, we note $A^*$ the group of its invertible elements.\\
For $E _2$ a finite extension of a local field $E_1$, we note respectively $Tr _{E _2 / E_1}$ and $N _{E _2 / E_1}$ the trace and norm of $E _2$ over $E _1$.\\
 We also note $Gal( E_2 / E_1)$ the Galois group of $E_2$ over $E_1$ when ${E _2 / E_1}$ is Galois, otherwise we note $Aut _{E _1} ({E _2} / {E _1})$ the group of automorphisms of the algebra $E _2$ over $E _1$. \\
 Moreover if $E _2$ is quadratic over $E _1$, we note $\eta _{E _2 / E_1}$ the nontrivial character of ${E _1}^*$ with kernel $N _{E _2 / E_1}({E _1}^*)$.\\
For $n$ a positive integer, we note ${GL _n (K)}_{+}$, the subgroup with index two of ${GL _n (K)}$, of matrices whose determinant is a norm of $K$ over $F$.\\
For $\Pi$ a representation of a group $G$, we note $\pi$ its class, and ${\Pi}^{\vee}$ its smooth contragedient when $\Pi$ is a smooth representation of a totally disconnected locally compact group.\\
If $\phi$ is an automorphism of $G$, we note ${\Pi}^{\phi}$ the representation of $G$ given by ${\Pi}\circ {\phi}$.\\
If $H$ is a subgroup of $G$, and $\mu$ is a character of $G$, we say that a representation $\Pi$ of $G$ is $\mu$-distinguished for ( with respect to) $H$ if there exists on the space of $ \Pi$ a linear functionnal $L$ verifying for $h$ in $H$, $L \circ \Pi (h) = \mu (h) L$.\\
If $ \mu$ is trivial, we say that $ \Pi$ is distinguished for $H$.\\

\subsection{ Quadratic extensions of $K$}

For $L$ a quadratic extension of $K$, three cases arise:

\begin{enumerate}

 \item  $L /F$ is biquadratic ( hence Galois), it contains $K$ and two other quadratic extensions $F$, $K'$ and $K''$. 

\begin{figure}[here]
 \centering
\includegraphics[width=2.1in,height=0.75in]{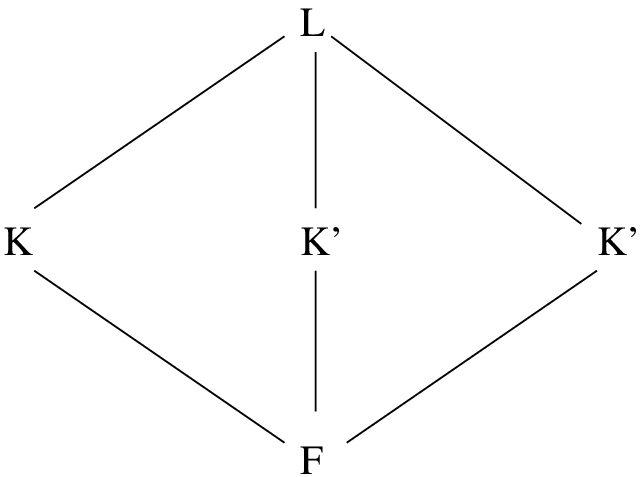}
\caption{}
\label{fig1}
\end{figure}

Its Galois group is isomorphic with $\mathbb{Z} / 2 \mathbb{Z} \times \mathbb{Z} / 2 \mathbb{Z}$, its non trivial elements are conjugations of $L$ over $K$, $K'$ and $K''$.\\

The conjugation $L$ over $K$ extend those of $K'$ and $K''$ over $F$.

 \item $L/F$ is cyclic with Galois group isomorphic with $\mathbb{Z} / 4 \mathbb{Z}$.

 \item $L/F$ non Galois. Then its Galois Closure $M$ is quadratic over $L$ and the Galois group of $M$ over $F$ is dihedral with order 8.\\
 To see this, we consider a morphism $\tilde{\theta}$ from $L$ to $\bar{F}$ which extends $\theta$. Then if $L' = \tilde{\theta}(L)$, $L$ and $L'$ are distinct, quadratic over $K$ and generate $M$ biquadratic over $K$. $M$ is the Galois closure of $L$ because any morphism from $L$ into $\bar{F}$, either extends $\theta$, or the identity map of $K$, sothat its image is either $L$ or $L'$, so it is always included in $M$.\\
Finally the Galois group $M$ over $F$ cannot be abelian ( for $L$ is not Galois), it is of order $8$, and it's not the quaternion group which only has one element of order $2$, whereas here the conjugations of $M$ over $L$ and $L'$ are of order $2$. Hence it is the dihedral group of order 8.\\ 
We deduce from this the folowing lattice:
 
\begin{figure}[h!]
 \centering
\includegraphics[width=3in,height=1.2in]{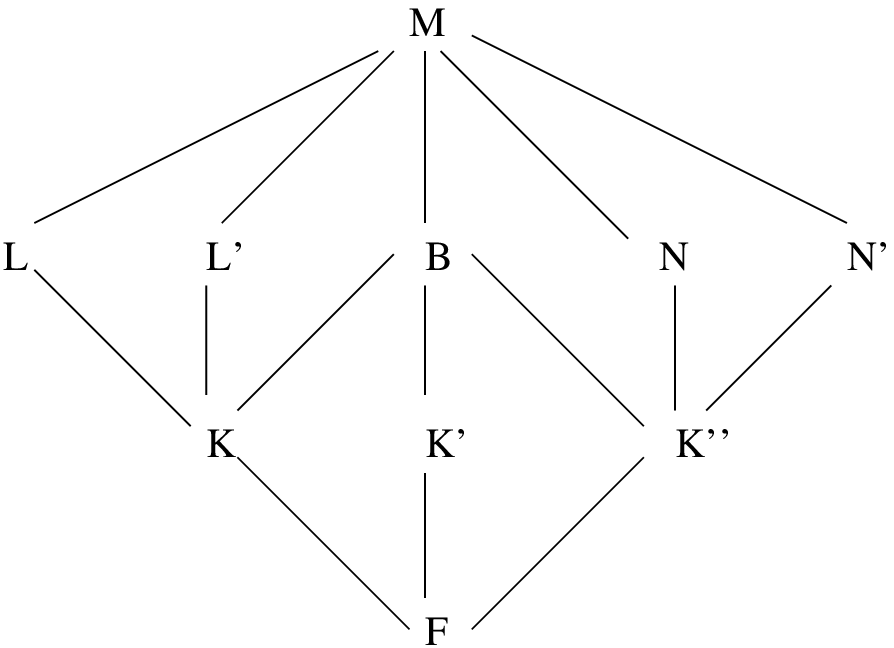}
\caption{}
\label{fig2}
\end{figure}

Here $M/K'$ is cyclic of degree 4, $M/K$ and $B/F$ are biquadratics.

\end{enumerate}

In the case $p$ odd, $F$ has exactly three quadratic extensions which generate its unique biquadratic extension.
If there exists $L$ non Galois over $F$, then it implies that the cardinal $q$ of the residual  field $F$ verifies $q\equiv 3[4]$, and $M$ is generated over $L$ by a primitive fourth root of unity in $\bar{F}$.

\subsection{Quadratic characters}  

We wish to calculate how $\eta _{L/K}$ restricts to $F^*$ in the following two cases. 

\begin{enumerate}

 \item If $L$ is biquadratic over $F$, then $\eta _{L /K}$ has a trivial restriction to $F^*$.\par
Indeed, we have $N _{{L}/K} ({K'}^*) = N _{K'/F}({K'}^*)$ and $N _{{L}/K} ({K''}^*) = N _{K''/F} ({K''}^*)$ because the conjugation of $L$ over $K$ extend those of $K'$ and $K''$ over $F$.\par
 Both these groups are distict from local class field theory and of index 2 in $F^*$, sothat they generate this latter, but both are contained in $N _{{L}/K} ({L}^*)$ which therefore contains $ F ^*$.\par
In other words $\eta _{L /K}$ restricts trivially to $F^*$.

\item If $L$ is cyclic over $F$, then ${\eta _{L /K}} _{|F^*}$ is non trivial.\par
If it wasn't the case, $F^*$ would be contained in $N _{L/K} (L^*)$, and composing with $N _{K/F}$ on both sides, ${F^*}^2$ would be a subgroup of $N _{L/F} (L^*)$. But ${F^*}^2$ and $N _{L/F} (L^*)$ have both index $4$ in $F^*$ and give different quotients ( $\mathbb{Z} / 2 \mathbb{Z} \times \mathbb{Z} / 2 \mathbb{Z}$ for the first and $\mathbb{Z} / 4 \mathbb{Z}$ for the second), sothat one cannot be contained in the other.

\end{enumerate}

\subsection{Weil's representation}

Let $L$ be a quadratic extension of $K$, then for any character $\omega$ of $L^*$ , we associate an irreducible representation $\Pi (\om)$ of $GL_2(K)$ ( cf.\cite{JL}), with central character $\om _{|K^*} \eta _{L/K}$.\\ 
If $\omega$ is regular for $N _{L/K}$, then $\Pi (\om)$ is supercuspidal, otherwise there exists a character $\mu$ of $K^*$ sothat $\om = \mu \circ N _{L/K}$, and then $\pi (\om)$ is the the principal series $\pi (\mu , \mu \eta _{L/K})$.\\

 The conjugation of $K$ over $F$ naturally extends to an involutive automorphism of $GL_2(K)$ which we also note $\theta$.\\
 
Here we want to determinate ${\Pi (\om)}^{\theta}$.\\ 
Suppose there exists $\tilde{\theta}$ an element of $Aut _{F}(L/F)$ which extends $\theta$, we then have:

\begin{prop}
If $\theta$ extends to an element $\tilde{\theta}$ of $Aut _{F}(L/F)$, then ${\Pi (\om)}^{\theta}$ is isomorphic with $\Pi ({\om}^{\tilde{\theta}})$ .

\end{prop}

\textbf{Proof}:\\
 Following \cite{JL}, $\Pi (\om)$ is the induced of $ r(\om,\psi _F)$ from $Gl_2 (K) _{+}$ the group $GL_2 (K)$.\par
 The space of this representation $ \textbf{S}( L, \om)$ is constitued by the continuous functions $f$ with compact support from $L^*$ to complex numbers verifying, for $x$ in $L$ and $y$ in $Ker(N _{L/K})$, $ f (xy) = {\om }^{-1} (y)f(x) $.\par  
 Then the mapping which associates to $f$ in $ \textbf{S}( L, \om)$ the function $f^{\theta}= f \circ \tilde{\theta}$ is an equivariant morphism between $ r(\om,\psi _F)^{\theta}$ and $ r(\om ^{\tilde{\theta}},\psi _F)$.\par
We then see ${\pi (\om)}^{\theta} = [Ind^{GL_2 (K)}_{GL_2  (K) _{+}} (r(\om,\psi _F))]^{\theta} \approx [Ind^{GL_2 (K)}_{GL_2 (K) _{+}} (r(\om,\psi _F) ^{\theta})] \approx Ind^{GL_2 (K)}_{GL_2 (K) _{+}} (r(\om ^{\tilde{\theta}},\psi _F)) = \pi (\om ^{\tilde{\theta}})$ where $Ind^{GL_2 (K)}_{GL_2 (K) _{+}}$ designs the induced representation from  $GL_2 (K) _{+}$ to $GL_2 (K)$.\\
\par

\textbf{Remark : }\\
It is not always true that $\theta$ extends to an element $\tilde{\theta}$ of $Aut _{F}(L/F)$.\\
For instance, take $F$ local with residual characteristic $q \equiv 3[4]$, and let $\pi _{F}$ be a prime element ( generating the maximal ideal of the integers ring).\\
We choose $K = F( {\pi _{F}}^{1/2})$ and $L = F( {\pi _{F}}^{1/4})$.\\
 Let $\tilde{\theta}$ be a $F$ linear morphism extending $\theta$ to $L$, with values in $\bar{F}$. Then $\tilde{\theta} ({\pi _{F}}^{1/4})= i {\pi _{F}}^{1/4}$, where $i$ is a primitive fourth root of unity.\\
Indeed $i^2=-1$, because   $\tilde{\theta} ({\pi _{F}}^{1/2})= \theta ({\pi _{F}}^{1/2}) = -{\pi _{F}}^{1/2}$.\\
Moreover, $i$ cannot be in $L$: indeed, this element is a root of unity with order prime to $q$, thus it would implie that the residual field of $L$, which is the one of $F$ as $L/F$ is totally ramified, contains a primitive fourth root of unity.\\
This cannot happen because $4$ does not divide $q-1$.\\ 
We conclude that any $F$ linear morphism extending $\theta$ to $L$, sends $L$ onto $ F( i {\pi _{F}}^{1/4})$ which is distinct from $L$ , and hence cannot be in $Aut _{F}(L/F)$.

\section{ Representations distinguished by a character}

\subsection{ Definitions and preliminary results }

The following theorem due to Y.Flicker \cite{AT} ( th. 1.3) will be of constant use.

\begin{thm}

Let $\Pi$ be an irreducible admissible representation of $GL_2 (K)$, such that $c_ \pi$ is trivial on $F^*$. Then ${\pi}^{\theta} ={\pi}^{\vee}$ if and only if $\pi$ is distingushed or $\eta _{K/F}$-distinguished for $GL_2 (F)$.
 
\end{thm}

Let $GL_2 (F) _{+}$ be the subgroup of index two in $GL_2 (F)$, it is clear that if a representation of $GL_2 (K)$ is distinguished or $\eta _{K/F}$-distinguished for $GL_2 (F)$, it is distinguished for $GL_2 (F) _{+}$, the reverse is true. 

\begin{prop}(cf.\cite{P}, p.71)

A representation of $GL_2 (K)$ is distinguished or $\eta _{K/F}$-distinguished for $GL_2 (F)$ if and only if it is distinguished for $GL_2 (F) _{+}$.

\end{prop}

\textbf{Proof}:\\
We show the non trivial implication.\par
Let $s$ be an element of $GL_2 (F)$ whose determinant is not a norm and let $\Pi$ be a $GL_2 (F) _{+}$-distingushed representation.\par
Let $L _{+}$ be the $GL_2 (F) _{+}$-invariant linear form on the space of $\Pi$, two cases arise:\

\begin{enumerate}
 
\item  If ${\Pi}^{\vee}  (s) L _{+}= - L _{+}$, then for $h$ in $GL_2 (F)\backslash GL_2 (F) _{+}$, we have ${\Pi}^{\vee} (h) L _{+} = {\Pi}^{\vee} (h) {\Pi}^{\vee} (s^2) L _{+}  = {\Pi}^{\vee} (hs) {\Pi}^{\vee} (s) L _{+}  = -L _{+}$ because $h s$ is in $GL_2 (F) _{+}$ ( here we also note ${\Pi}^{\vee}$ the non smooth contragedient). ${\Pi}$ is therefore  $\eta _{K/F}$-distinguished for $GL_2 (F)$.

\item Otherwise ${\Pi}^{\vee}  (s) L _{+} \neq - L _{+}$, and $L _{+} + {\Pi}^{\vee}  (s) L _{+}$ is fixed under the action of $GL_2 (F)$.

\end{enumerate}

Theorem 3.1 takes the following form:

\begin{thm}

Let $\Pi$ be an irreducible admissible representation of $GL_2 (K)$. Then $\Pi$ is $GL_2 (F) _{+}$-distinguished if and only if ${\pi}^{\theta} = {\pi}^{\vee}$ and $c_ \Pi$ restricts trivially to $F^{*}$. 
 
\end{thm}

\subsection{Description of the $GL_2 (F) _{+}$-distinguished representations}

\begin{thm} 

 A supercuspidal dihedral representation $\Pi$ of $GL_2 (K)$ is $GL_2 (F) _{+}$-distinguished if and only if there exists a quadratic extension $L$ of $K$ biquadratic on $F$, and a multiplicative character $\omega$ of $L$ trivial on $N _{L/K'}({K'}^*)$ or on $ N _{L/K''}({K''}^*)$, such that $\pi = \pi (\omega)$.

\end{thm}

\textbf{Proof:}\par

Let $L$ be a quadratic extension of $K$ and $\omega$ a regular multiplicative of $L$ such that $\pi = \pi (\omega)$, we note $\sigma$ the conjugation of $L$ over $K$, three cases show up:\

\begin{enumerate}

\item $L/F$ si biquadratic.\\
we note $\sigma '$ the conjugation of $L$ over $K'$ and $\sigma ''$ the conjugation of $L$ over $K''$, $\sigma '$ and $\sigma ''$ both extend $\theta$, and thus can play $\tilde{\theta}$'s role in proposition 1.1.\\
 The condition ${\pi}^{\vee} = {\pi}^{\theta}$ which one can also read $\pi ({\omega}^{-1})= \pi ({\omega}^{\tilde{\theta}})$, is then equivalent to ${\omega}^{\sigma '}= {\omega}^{-1}$ or ${\omega}^{\sigma ''}={\omega}^{-1}$.\\
 This is equivalent to $\om$ trivial on $N _{L/K'}({K'}^*)$ and on $ N _{L/K''}({K''}^*)$.\\
 As $\eta _{L/K}$, $\eta _{L/K'}$, and $\eta _{L/K''}$ are trivial  $F^*$, we have ${c _{\pi}} _{ | F^*}= {\omega} _{| F^*} {\eta _{L/K}} _{| F^*} = 1$ for such a representation .

\item $L/F$ is cyclic, the regularity of $\omega$ makes the condition $\pi ({\omega}^{-1}) = {\pi ({\omega}^{\tilde{\theta}})}$ impossible.\\
Indeed one would have ${\omega}^{\tilde{\theta}}= {\omega}^{-1}$, which would implie ${\omega}^{\sigma}= \omega$ for ${\sigma}^2 = \theta$, and so $\omega$ would be trivial on the kernel of $N _{L/K}$ from Hilbert's theorem 90.\\
$\Pi$ can therefore not be $GL_2 (F) _{+}$-distniguished.\

\item $L/K$ is not Galois ( which implies $q \equiv 3[4]$ in the case p odd), we note again $\theta$ the conjugation of $B$ over $K'$ which extends the one of $K$ over $F$.\\
 Let ${\pi}_{B/K}$ be the representation of ${GL _2 (B)}$ which is the base change lifting of ${\pi}$ to $B$. As ${\pi}_{B/K}= \pi ({\omega} \circ N _{M/L})$, if $ {\omega} \circ N _{M/L} = {\mu}\circ N _{M/B} $ for a character $\mu$ of $B^*$, then $\pi(\omega) = \pi(\mu)$ ( cf.\cite{GL}, (3) p.123) and we are brought back to case 1.\\
 Otherwise $ {\omega} \circ N _{M/L}$ is regular for $ N _{M/B} $.\\
 If ${\pi}$ was ${GL _2 (F)}_{+}$-distinguished, taht is ${\pi}^{\theta} ={\pi}^{\vee}$ and ${c _{\pi}} _{ | F^*}$, we would have ${\pi}_{B/K}^{\theta}= {\pi}_{B/K}^{\vee}$ and $c _{{\pi}_{B/K}}= c_{\pi} \circ N _{B/K}$ from theorem 1 of \cite{GL}.\\
As $N _{B/K} ({K'} ^*)= N _{K'/F} ({K'} ^*)$ for the conjugation of $B$ over $K$ extends that of $K'$ over $F$, one would deduce that $c _{{\pi}_{B/K}}$ would be trivial on ${K'} ^*$ and theorem 2.2 would implie that ${\pi}_{B/K}$ would be ${GL _2 (K')}_{+}$-distinguished.\\
 That would contradict case 2 because $M/K'$ is cyclic.

\end{enumerate}

\section{Distingushed representations}

We described in the previous section the supercuspidal dihedral representations of $GL_2(K)$ which are ${GL _2 (F)}_{+}$ distinguished.\\  
We want to characterize those who are ${GL _2 (F)}$-distinguished among them.\\

\subsection {Definitions and useful results}

We refer to \cite{JL} for definitions and basic properties of $\epsilon$ factors attached to an irreducible admissible representation of $GL_2(K)$, and to \cite{T} for those of $\epsilon$ factors attached to a multplicative character of a local field.\\
 The $\epsilon$ used here for representations of $GL_2(K)$ is the one described in \cite{JL} evaluated at $s = 1/2$ and the $\epsilon$ attached to a multplicative character of a local field is Langlands'$\epsilon _L$ described in \cite{T}.\\

We will use the three following results.\\ \par

The first, due to J.Hakim can be found in \cite{H}, page 8.\\
 Here we repaced $\gamma$ with $\epsilon$ because both are equal for supercuspidal representations:

\begin{thm}

Let $\Pi$ be a supercuspidal irreducible representation of $GL_2(K)$, and $\psi$ a nontrivial character of $K$ trivial on $F$. Then  $\Pi$ is distinguished if and only if $\epsilon (\Pi \otimes \x , \psi)=1 $ for every character $\x$ of $K^*$ trivial on $F^*$. 

\end{thm}

The second, due to Fröhlich and Queyrut, is in \cite{FQ}, page 130 :

\begin{thm}
 Let $L_2$ be a quadratic extension of $L_1$ which is a quadratic extension of $\mathbb{Q} _p$ , then if $\psi _{L_2}$ is the standard character of $L_2$ and if $\Delta$ is an element of ${L_2}^*$ with $Tr _{L_2 /L_1} (\Delta) = 0 $, we then have $\epsilon ( \x , \psi_ {L_2} )= \x (\Delta) $ for every character $\x$ of $L_2 ^*$ trivial on $L_1 ^*$ .  

\end{thm}

The third is a corollary of proposition 3.1 of \cite{AT}:

\begin{prop}
There exists no supercuspidal representation of $GL_2(K)$ which is distinguished and $\eta _{K/F}$-distinguished at the same time.

\end{prop}

\section{Description of distingushed representations}

\begin{thm} 

 A dihedral supercuspidal representation $\Pi$ of $GL_2 (K)$ is $GL_2 (F)$-distinguished if and only if there exists a quadratic extension $L$ of $K$ biquadratic over $F$, and a regular multiplicative character  $\omega$ of $L$ trivial on ${K'}^*$ or on ${K''}^*$, such that $\pi = \pi (\omega)$.

\end{thm}

\textbf{Proof:}\par

From the second section, we can suppose that $\pi = \pi (\om)$, for $\omega$ a regular multiplicative character of a quadratic extension $L$ of $K$ biquadratic over $F$, with $\om$ trivial on $N _{L/K'}({K'}^*)$ or on $ N _{L/K''}({K''}^*)$.\\  
\par 

Let $\psi _K$ be the standard character of $K$, $\psi _L$ the one of $L$, and $a$ a non null element of $K$ such that $Tr _{K/F} (a)=0$, which implies $Tr _{L/K'} (a)= Tr _{L/K''} (a)=0$.\\

we note $(\psi _K)_a$ the character trivial on $F$ given by $(\psi _K)_a (x) = (\psi _K)(ax)$.\\

To see if $\pi(\omega)$ is distinguished, we use Hakim's criterion ( th.3.1).\\

So let $\chi$  be a character of $K^*$ trivial on $F^*$, we have $ \pi(\omega) \otimes \chi = \pi(\omega \times \chi \circ N _{L / K} )$ and we note $ \mu =  \omega \times \chi \circ N _{L/ K}$.\\
\\
 
\textbf{i)} if ${\omega} _{|{K'}^*}=1$: on a $\epsilon (\pi(\omega) \otimes \chi , (\psi _K)_a)= \epsilon ( \pi (\mu) , ({\psi} _{K}) _a)= \epsilon (\pi (\mu) , \psi _K) \mu(a) {\eta} _{L/K} (a)$. Now $\epsilon (\pi (\mu) , \psi _K)= \lambda (L/K, \psi _K) \epsilon ( \mu , \psi _L)$ ( cf.\cite{JL} p.150), where the Langlands-Deligne factor $\lambda (L/K, \psi _K)$ equals $\epsilon ( \eta _{L/K} , \psi _K)$ divided by its module.\\

 As ${\eta _{L/K}} _{|F^*}=1$ et $\mu _{|{K'}^*}=1$, from theorem 4.2, we have that $\epsilon ( \mu , \psi _L)= \mu (a)$ and $\epsilon ( \eta _{L/K} , \psi _K) = \eta _{L/K} (a)$.\\ 

We deduce that $\epsilon (\pi(\omega) \otimes \chi , (\psi _K)_a)= {\mu(a)}^2 {{\eta} _{L/K} (a)}^2 = 1$ for $a^2$ is in $F$. $\pi(\omega)$ is therefore distinguished.\\
\\

\textbf{ii)} If ${\omega} _{|{K'}^*}= \eta _{L/K'}$ : Let $\chi '$ be a character of $K^*$ which extends ${\eta} _{K/F}$, then $\chi ' \circ N _{L/K}$ equals $\eta _{K/F} \circ N _{K'/F}$ on $K'$ because the conjugation of $L$ over $K$ extends the one of $K'$ over $F$.\\
But $\eta _{K/F} \circ N _{K'/F}$ is trivial on the image of $N _{L/K'}$ from the identity $N _{K'/F} \circ N _{L/K'} = N _{K/F} \circ N _{L/K} $, but not trivial for $N _{K' /F}$ is not the kernel $N _{K /F}(K^*)$ of $\eta _{K/F}$ from local class field theory.\\
Thus $\omega \times \chi ' \circ N _{L/K}$ is trivial on $K'$, and we deduce that $\pi (\omega)\otimes \chi ' = \pi (\omega \times \chi ' \circ N _{L/K})$ is distinguished from i).\\
 This implies that $\pi (\omega)$ is ${\eta} _{K/F}$-distinguished and thus not distinguished from proposition 3.1.

The cases  ${\omega} _{|{K''}^*}=1$ and ${\omega} _{|{K''}^*}=\eta _{L/K''}$ are handled as well.\\

\section{The principal series}

Representations of the principal series of $GL_2 (K)$ distinguished for $GL_2 (F)$ are well known, and described for example in proposition 4.2 of \cite{AT}.\par

 The result is the following:

\begin{prop}

Let $\lambda$ and $\mu$ be two characters of $K^*$, whose quotient is not the module of $K$ or its inverse. The principal series representation $\Pi (\lambda, \mu)$ of $GL_2 (K)$ is $GL_2 (F)$-distinguished either when $\lambda= {\mu}^{- \theta}$ or when $\lambda$ and $\mu$ have a trivial restriction to $F^*$.

\end{prop}

Now one can construct the principal series $\Pi (\lambda, \mu)$ via the Weil representation ( cf.\cite{B} p.523 à 557), in this case $(\lambda, \mu)$ identifies with a character of $K^* \times K^*$.\par
This way of parametrising irreducible representations of $GL_2 (K)$ with multiplicatice characters of two-dimensional semi-simple commutative algebras over $K$, includes the principal series (for the algebra $K \times K$) and the dihedral representations (for quadratic extensions of $K$).\\
\par
Let $L$ be a quadratic extension of $K$ biquadratic over $F$, as we are here interested with $GL_2 (F)$-distinguishedness, we consider the following $F$-algebras.\par

\begin{enumerate}

 \item  \textbf{the algebra $K \times K$}\par

One note $Aut _F (K \times K)$ its automorphisms group. The elements of this group are $(x,y)\mapsto(x,y)$, $(x,y)\mapsto(y,x)$, $(x,y)\mapsto(x^{\theta}, y^{\theta})$, $(x,y)\mapsto(y^{\theta}, x^{\theta})$, and $Aut _F (K \times K)$ is isomorphic with $\mathbb{Z} / 2 \mathbb{Z} \times \mathbb{Z} / 2 \mathbb{Z}$.\par
The three sub-algebras fixed by non trivial elements of $Aut _F (K \times K)$ are $K$ via the natural diagonal inclusion, the twisted form ${\tilde{K}}$ of $K$ given by $x\mapsto(x,x^{\theta})$, and $F\times F$.

\item \textbf{l'algèbre $L$}\par

The group $Gal(B _{F} /K )$ of its automorphisms is isomorphic with $\mathbb{Z} / 2 \mathbb{Z} \times \mathbb{Z} / 2 \mathbb{Z}$.\par

The three sub-algebras fixed by non trivial elements of $Gal(B _{F} /K )$ are $K$, $K'$ et $K''$.

\end{enumerate}

We then observe that proposition 5.2 for the principal series has the same statement that the one for theorem 4.1:

\begin{prop}

  A principal series representation $\Pi(\lambda, \mu)$  of $GL_2 (K)$ is $GL_2 (F)$-distinguished if and only if the multiplicative character $(\lambda, \mu)$ is trivial one the invertible elements of one of the two intermediate sub-algebras of $K \times K$ distinct from $K$.

\end{prop}

We now study dihedral non supercuspidal representations.\\
Let $\Pi$ be such a representation, there exists a quadratic extension $L$ over $K$ and a non regular multiplicative character $\omega$ of $L$ such that $\pi= \pi (\om)$.\\
 If $\mu $ is a character of $K^*$ such that $\om = \mu \circ N _{L/K}$, then $\pi= \pi( \mu , \mu \eta _{L / K} )$.\\ 

Three cases arise:

\begin{enumerate}
 
\item For $L$ biquadratic over $F$, we show that $\omega$ restricts trivially to ${K'}^*$ or ${K''}^*$ if and only if $( \mu , \mu \eta _{L / K} )$ restricts trivially to ${\tilde{K}}^*$ or to ${F^* \times F^*}$.

We have the following equivalences:
\begin{itemize}

 \item $\omega({K'}^*)=1 \Leftrightarrow \mu(N _{L /K} ({K'}^*))=1 \Leftrightarrow \mu(N _{K' /F} ({K'}^*))=1 $ because the conjugation of $L /K$ extends the one of $K' /F$, and so $\omega({K'}^*)=1 \Leftrightarrow$ $\mu_{|F ^*}=1$ or $\eta _{K' /F}$.

\item $\omega({K''}^*)=1 \Leftrightarrow \mu(N _{L /K} ({K''}^*))=1 \Leftrightarrow \mu(N _{K'' /F} ({K''}^*))=1 $ because the conjugation of $L /K$ extends the one of $K'' /F$, and so $\omega({K''}^*)=1 \Leftrightarrow$ $\mu_{|F ^*}=1$ or $\eta _{K'' /F}$.

\item $( \mu , \mu \eta _{L / K} )$ trivial on ${\tilde{K}}^*$ $\Longleftrightarrow {\mu}^{\theta} \mu \eta _{L / K}=1 \Longleftrightarrow \mu \circ N _{K /F} = \eta _{L / K}$.\par
 We deduce that $\mu \circ N _{L /F} = 1$, which implies that $\mu_{|F ^*}=1$, $\eta _{K /F}$, $\eta _{K' /F}$ or $\eta _{K'' /F}$, but the first two possibilities are excluded by the identity $\mu \circ N _{K /F} = \eta _{L / K}$.\par
Conversely if $\mu_{|F ^*}=\eta _{K' /F}$ or $\eta _{K'' /F}$, then $\mu \circ N _{K /F}$ is a character of order two of $K^*$ which cannot be trivial from local class field theory. As the equalities $N _{L / F}= N _{L / K'} \circ N _{K'/ K}= N _{L / K''} \circ N _{K''/ K}$ implie that $\mu \circ N _{K /F}$ is trivial on $N _{L / K} ({L} ^*)$, it is therefore $\eta _{L/ K}$.\par
Eventually $( \mu , \mu \eta _{L / K} )$ trivial on ${\tilde{K}}^*$ $\Longleftrightarrow$ $\mu_{|F ^*}=\eta _{K' /F}$ or $\eta _{K'' /F}$.

\item Also $( \mu , \mu \eta _{L / K} )$ trivial on ${F^* \times F^*}$ $\Longleftrightarrow$ $\mu_{|F ^*}=1$ because we have already seen that $\eta _{L / K} $ is trivial on $F ^*$.

\item Finally these equivalences show that $\omega({K'}^*)=1$ or $\omega({K''}^*)=1$ $\Longleftrightarrow$ $( \mu , \mu \eta _{L / K} )$ trivial on ${\tilde{K}}^*$ or on ${F^* \times F^*}$.

\end{itemize}

\item If $L$ is cyclic over $F$. One shows that $\pi(\om)= \pi( \mu , \mu \eta _{L / K} )$ is distinguished if and only if  $\mu _{|F^*}$ generates the cyclic group of the characters of $F^*/ N _{L /F}(L^*)$.

\begin{itemize}

\item It is not possible for $( \mu , \mu \eta _{L / K} )$ to be trivial on ${F^* \times F^*}$ because we saw in the preliminaries that $\eta _{L / K} $ is not trivial on $F ^*$.

\item $( \mu , \mu \eta _{L / K} )$ trivial on ${\tilde{K}}^*$ $\Longleftrightarrow \mu \circ N _{K /F} = \eta _{L / K}$.\par
We deduce as before that $\mu$ is a character of $F^* / N _{L /F}(L^*)$.\\
As $F^* / N _{L /F}(L^*)$ is cyclic of order four, the same is true for its characters group.\\
 As $\mu \circ N _{K /F} = \eta _{L / K}$, we deduce that $\mu _{|F^*}$ is non trivial, moreover if $\mu _{|F^*}$ was of order 2, it would be equal to $\eta _{K / F}$ which is the unique element with order 2 of the characters group of $F^* / N _{L /F}(L^*)$, which contradicts $\mu \circ N _{K /F} = \eta _{L / K}$.\\
 We deduce that $\mu \circ N _{K /F} = \eta _{L / K}$ $\Longrightarrow$ $\mu _{|F^*}$ generates the dual group of $F^* / N _{L /F}(L^*)$. 

\item Conversely if $\mu _{|F^*}$ is of order four in the dual group of $F^* / N _{L /F}(L^*)$, we deduce that $\mu \circ N _{K /F}$ is a character of $K^*$ trivial on $N _{L /K}(L^*)$, but not trivial because it would implie $\mu _{|F^*}=1$ or $\eta _{K/F}$, namely $\mu _{|F^*}$ with order less than 2.\\
On conclude that $\mu \circ N _{K /F} = \eta _{L / K}$.

\end{itemize}

\item If $L$ is not Galois over $F$. As $\omega \circ N _{M/K} = \mu \circ N _{M/F} = \mu ' \circ N _{M/B}$, where $\mu '= \mu \circ N _{B/F}$, we conclude as in the case 3. of the proof of theorem 2.3 that $\pi (\om) =  \pi (\mu ')$ and we are brought back to case 1. because $B$ is biquadratic over $F$. 

\end{enumerate}

Thus we have the following general theorem:

\begin{thm}

A dihedral representation $\Pi$ of $GL_2 (K)$ is $GL_2 (F)$-distinguished if and only if $\pi = \pi (\omega)$ for some multiplicative character $\omega$ of a quadratic extension $L$ over $K$ verifying i) or ii):

i) $L/F$ is biquadratic , and $\omega _{|{K'}^*}=1$ or $\omega _{|{K''}^*}=1$,

ii) $L/F$ is cyclic and  $\omega = \mu \circ N _{L/K}$ for $\mu$ a character of $K^*$ whose restriction to $F^*$ generates the dual group of $F^* / N _{L/F} (L^*)$.

\end{thm}


\begin{thebibliography}{9}


\bibitem[A-T]{AT} U.K.Anandavardhanan and R.Tandon , \textit{ On Distinguishness }, Pac. J. Math., \textbf{206} (2002), 269-286.

\bibitem[B]{B} Daniel Bump, \textit{ Automorphic Forms and Representations }, Cambridge studies in advanced Mathematics 55, Cambridge University Press,1998.


\bibitem[F-Q]{FQ} A.Frohlich and J.Queyrut, \textit{ On the functional equation of the Artin L-Function for characters of real representations}, Inv. Math., \textbf{20} (1973), 125-138. 


\bibitem[G-L]{GL} P.Gérardin and J-P.Labesse, \textit{ The solution of a base change problem for $GL(2)$ }, in Automorphic Forms, Representations and $L-$functions (Corvallis), AMS Proc. Symp. Pure. Math., \textbf{33} (1979), 155-133.


\bibitem[H]{H} J.Hakim, \textit{ Distinguished p-adic representations}, Duke Math. J., \textbf{62} (1991), 1-22.


\bibitem[J-L]{JL} H.Jacquet and R.Langlands , \textit{ Automorphic forms on GL(2) }, Lect. Notes Math., \textbf{114}, Springer, 1970.

\bibitem[P]{P} M.N. Panichi, \textit{Caractérisations du spectre tempéré de  $GL _n (\mathbb{C}) / GL _n (\mathbb{R})$}, Thèse de Doctorat, Université Paris 7, 2001.


\bibitem[Se91]{S1} J-P.Serre, \textit{Petit cours d'arithmétique }, PUF, 1991.


\bibitem[Se98]{S2} J-P.Serre, \textit{Représentations linéaires des groupes finis}, Hermann, 1998, 74-76.


\bibitem[T]{T} J.Tate, \textit{ Number theoritic background }, in Automorphic Forms, Representations and $L-$functions (Corvallis), AMS Proc. Symp. Pure. Math., \textbf{33} (1979), 3-26.

\bibitem[W]{W} A.Weil, \textit{ Basic Number Theory }, Springer-Verlag, 1973.



\end{thebibliography}
\end{document}